\documentclass[12pt,a4paper,reqno]{amsart}
\usepackage{amsmath}
\usepackage{amsfonts}
\usepackage{amssymb}
\numberwithin{equation}{section}

\theoremstyle{plain}

\newtheorem{theorem}[subsection]{Theorem}
\newtheorem{proposition}[subsection]{Proposition}

\newtheorem{conjecture}[subsection]{Conjecture}

\theoremstyle{definition}

\newtheorem{definition}[subsection]{Definition}

\newtheorem{question}[subsection]{Question}
\newtheorem{dichotomy}[subsection]{Dichotomy}

\renewcommand{\leq}{\leqslant}
\renewcommand{\geq}{\geqslant}

\newsavebox{\proofbox}
\savebox{\proofbox}{\begin{picture}(7,7)%
  \put(0,0){\framebox(7,7){}}\end{picture}}
\def\boxeq{\tag*{\usebox{\proofbox}}}

\newcommand{\md}[1]{\ensuremath{(\mbox{mod}\, #1)}}

\def\endproof{\hfill{\usebox{\proofbox}}}
\def\E{\mathbb{E}}
\def\Z{\mathbb{Z}}

\def\P{\mathbb{P}}

\def\ni{\noindent}
\def\vs{\vspace{11pt}}

\begin{document}

\title{Long arithmetic progressions of primes 
}

\author{Ben Green}
\address{School of Mathematics, University Walk, Bristol BS8 1TW, England.
}
\email{b.j.green@bristol.ac.uk}

\begin{abstract}
This is an article for a general mathematical audience on the author's work, joint with Terence Tao, establishing that there are arbitrarily long arithmetic progressions of primes. 
\end{abstract}

\maketitle

\section{introduction and history}

\ni This is a description of recent work of the author and Terence Tao \cite{green-tao-primes} on primes in arithmetic progression. It is based on seminars given for a general mathematical audience in a variety of institutions in the UK, France, the Czech Republic, Canada and the US.\vs

\ni Perhaps curiously, the order of presentation is much closer to the order in which we discovered the various ingredients of the argument than it is to the layout in \cite{green-tao-primes}. We hope that both expert and lay readers might benefit from contrasting this account with \cite{green-tao-primes} as well as the expository accounts by Kra \cite{Kra} and Tao \cite{taosurvey1,taosurvey2}.\vs

\ni As we remarked, this article is based on lectures given to a general audience. It was often necessary, when giving these lectures, to say things which were not strictly speaking true for the sake of clarity of exposition. We have retained this style here. However, it being undesirable to commit false statements to print, we have added numerous footnotes alerting readers to points where we have oversimplified, and directing them to places in the literature where fully rigorous arguments can be found.\vs

\ni Our result is:
\begin{theorem}[G.--Tao]\label{mainthm1}
The primes contain arbitrarily long arithmetic progressions.\endproof
\end{theorem}
\ni Let us start by explaining that the truth of this statement is not in the least surprising. For a start, it is rather easy to write down a progression of five primes (for example $5,11,17,23,29$), and in 2004 Frind, Jobling and Underwood produced the example
\[ 56211383760397 + 44546738095860 k; \quad k = 0, 1, \ldots, 22.\]
of 23 primes in arithmetic progression. A very crude heuristic model for the primes may be developed based on the prime number theorem, which states that $\pi(N)$, the number of primes less than or equal to $N$, is asymptotic to $N/\log N$. We may alternatively express this as
\[ \P\big( \mbox{$x$ is prime} \; | \; 1 \leq x \leq N\big) \sim 1/\log N.\]
Consider now the collection of all arithmetic progressions
\[ x, x+ d, \dots, x + (k-1)d\]
with $x,d \in \{1,\dots,N\}$. Select $x$ and $d$ at random from amongst the $N^2$ possible choices, and write $E_j$ for the event that $x + jd$ is prime, for $j = 0,1,\dots, k-1$. The prime number theorem tells us that 
\[ \P(E_j) \approx 1/\log N.\]
If the events $E_j$ were independent we should therefore have
\[ \P(\mbox{$x, x+ d, \dots, x + (k-1)d$ are all prime}) = \P \big( \bigwedge_{j=0}^{k-1} E_j \big) \approx 1/(\log N)^k.\] 
We might then conclude that
\[ \# \{x,d \in \{1,\dots,N\}: x, x+ d, \dots, x + (k-1)d \; \mbox{are all prime} \; \} \approx \frac{N^2}{(\log N)^k}.\]
For fixed $k$, and in fact for $k$ nearly as large as $2\log N/\log \log N$, this is an increasing function of $N$. This suggests that there are infinitely many $k$-term arithmetic progressions of primes for any fixed $k$, and thus arbitrarily long such progressions.\vs

\ni Of course, the assumption that the events $E_j$ are independent was totally unjustified. If $E_0, E_1$ and $E_2$ all hold then one may infer that $x$ is odd and $d$ is even, which increases the chance that $E_3$ also holds by a factor of two. There are, however, more sophisticated heuristic arguments available, which take account of the fact that the primes $> q$ fall only in those residue classes $a \md{q}$ with $a$ coprime to $q$. There are very general conjectures of Hardy-Littlewood which derive from such heuristics, and a special case of these conjectures applies to our problem. It turns out that the extremely na\"{\i}ve heuristic we gave above only misses the mark by a constant factor:

\begin{conjecture}[Hardy-Littlewood conjecture on $k$-term APs]\label{hlconj}
For each $k$ we have 
\[ \# \{x,d \in \{1,\dots,N\} : x, x+ d, \dots, x + (k-1)d \;\; \mbox{\textup{are all prime}} \; \} = \frac{\gamma_k N^2}{(\log N)^k}(1 + o(1)),\]
where 
\[ \gamma_k = \prod_p \alpha^{(k)}_p\] is a certain product of ``local densities'' which is rapidly convergent and positive.
\end{conjecture}
\ni We have
\[ \alpha_p^{(k)} = \left\{ \begin{array}{ll} \frac{1}{p} \left(\frac{p}{p-1}\right)^{k-1} & \mbox{if $p \leq k$} \\ \left(1 - \frac{k-1}{p}\right) \left( \frac{p}{p-1}\right)^{k-1} & \mbox{if $p \geq k$}.\end{array}\right.\]
In particular we compute\footnote{For a tabulation of values of $\gamma_k$, $3 \leq k \leq 20$, see \cite{gh}. As $k \rightarrow \infty$, $\log \gamma_k \sim k \log \log k$.}
\[ \gamma_3 = 2\prod_{p \geq 3} \big(1 - \frac{1}{(p-1)^2}\big) \approx 1.32032\]
and 
\[ \gamma_4 = \frac{9}{2}\prod_{p \geq 5} \big(1 - \frac{3p-1}{(p-1)^3}\big) \approx 2.85825.\]
\ni What we actually prove is a somewhat more precise version of Theorem \ref{mainthm1}, which gives a lower bound falling short of the Hardy-Littlewood conjecture by just a constant factor.

\begin{theorem}[G.--Tao]\label{gtthm} For each $k \geq 3$ there is a constant $\gamma'_k > 0$ such that
\[ \# \{x,d \in \{1,\dots,N\} : x, x+ d, \dots, x + (k-1)d \;\; \mbox{\textup{are all prime}} \; \} \geq \frac{\gamma'_k N^2}{(\log N)^k}\] for all $N > N_0(k)$.\endproof
\end{theorem}

\ni The value of $\gamma'_k$ we obtain is very small indeed, especially for large $k$.\vs 

\ni Let us conclude this introduction with a little history of the problem. Prior to our work, the conjecture of Hardy-Littlewood was known only in the case $k = 3$, a result due to Van der Corput \cite{van-der-corput} (see also \cite{chowla}) in 1939. For $k \geq 4$, even the existence of infinitely many $k$-term progressions of primes was not previously known. A result of Heath-Brown from 1981 \cite{heath-brown1} comes close to handling the case $k = 4$; he shows that there are infinitely many $4$-tuples $q_1 < q_2 < q_3 < q_4$ in arithmetic progression, where three of the $q_i$ are prime and the fourth is either prime or a product of two primes. This has been described as ``infinitely many $3\frac{1}{2}$-term arithmetic progressions of primes''.

\section{the relative szemer\'edi strategy}\label{sec2}

\ni A number of people have noted that \cite{green-tao-primes} manages to avoid using any deep facts about the primes. Indeed the only serious number-theoretical input is a zero-free region for $\zeta$ of ``classical type'', and this was known to Hadamard and de la Vall\'ee Poussin over 100 years ago. Even this is slightly more than absolutely necessary; one can get by with the information that $\zeta$ has an isolated pole at $1$ \cite{taonotes1}.\vs

\ni Our main advance, then, lies not in our understanding of the primes but rather in what we can say about \textit{arithmetic progressions}. Let us begin this section by telling a little of the story of the study of arithmetic progressions from the combinatorial point of view of Erd\H{o}s and Tur\'an \cite{erdos-turan}.\vs

\begin{definition}
Fix an integer $k \geq 3$. We define $r_k(N)$ to be the largest cardinality of a subset $A \subseteq \{1,\dots,N\}$ which does not contain $k$ distinct elements in arithmetic progression.
\end{definition}
\ni Erd\H{o}s and Tur\'an asked simply: what is $r_k(N)$? To this day our knowledge on this question is very unsatisfactory, and in particular we do not know the answer to
\begin{question}
Is it true that $r_k(N) < \pi(N)$ for $N > N_0(k)$?
\end{question}
\ni If this is so then the primes contain $k$-term arithmetic progressions on density grounds alone, irrespective of any additional structure that they might have. I do not know of anyone who seriously doubts the truth of this conjecture, and indeed all known lower bounds for $r_k(N)$ are much smaller than $\pi(N)$. The most famous such bound is Behrend's assertion \cite{behrend} that 
\[ r_3(N) \gg N e^{-c\sqrt{\log N}};\]
slightly superior lower bounds are known for $r_k(N)$, $k \geq 4$ (cf. \cite{laba-lacey,rankin}).\vs

\ni The question of Erd\H{o}s and Tur\'an became, and remains, rather notorious for its difficulty. It soon became clear that even seemingly modest bounds should be regarded as great achievements in combinatorics. The first really substantial advance was made by Klaus Roth, who proved 

\begin{theorem}[Roth, \cite{roth}] We have $r_3(N) \ll N(\log \log N)^{-1}$.\endproof
\end{theorem}
\ni The key feature of this bound is that $\log \log N$ tends to infinity with $N$, albeit slowly\footnote{cf. the well-known quotation ``$\log\log\log N$ has been proved to tend to infinity with $N$, but has never been observed to do so''.}. This means that if one fixes some small positive real number, such as 0.0001, and then takes a set $A \subseteq \{1,\dots,N\}$ containing at least $0.0001 N$ integers, then provided $N$ is sufficiently large this set $A$ will contain three distinct elements in arithmetic progression.\vs

\ni The generalisation of this statement to general $k$ remained unproven until  Szemer\'edi clarified the issue in 1969 for $k = 4$ and then in 1975 for general $k$. His result is one of the most celebrated in combinatorics.

\begin{theorem}[Szemer\'edi  \cite{szemeredi-4,szemeredi}]
We have $r_k(N) = o(N)$ for any fixed $k \geq 3$.\endproof
\end{theorem}

\ni Szemer\'edi's theorem is one of many in this branch of combinatorics for which the bounds, if they are ever worked out, are almost unimaginably weak. Although it is in principle possible to obtain an explicit function $\omega_k(N)$, tending to zero as $N \rightarrow \infty$, for which 
\[ r_k(N) \leq \omega_k(N)N,\] to my knowledge no-one has done so. Such a function would certainly be worse than $1/\log_* N$ (the number of times one must apply the log function to $N$ in order to get a number less than 2), and may even be slowly-growing compared to the inverse of the Ackermann function.\vs

\ni The next major advance in the subject was another proof of Szemer\'edi's theorem by Furstenberg \cite{furst}. Furstenberg used methods of ergodic theory, and his argument is relatively short and conceptual. The methods of Furstenberg have proved very amenable to generalisation. For example in \cite{bergelson-leibman} Bergelson and Leibman proved a version of Szemer\'edi's theorem in which arithmetic progressions are replaced by more general configurations $(x + p_1(d),\dots, x + p_k(d))$, where the $p_i$ are polynomials with $p_i(\mathbb{Z}) \subseteq \mathbb{Z}$ and $p_i(0) = 0$. A variety of multidimensional versions of the theorem are also known. A significant drawback\footnote{A discrete analogue of Furstenberg's argument has now been found by Tao \cite{tao-discrete-ergodic}. It does give an explicit function $\omega_k(N)$, but once again it tends to zero incredibly slowly.} of Furstenberg's approach is that it uses the axiom of choice, and so does not give \textit{any} explicit function $\omega_k(N)$. \vs

\ni Rather recently, Gowers \cite{gowers4,gowers} made a major breakthrough in giving the first ``sensible'' bounds for $r_k(N)$. 

\begin{theorem}[Gowers]
Let $k \geq 3$ be an integer. Then there is a constant $c_k > 0$ such that 
\begin{equation}\boxeq r_k(N) \ll N(\log \log N)^{-c_k}.\end{equation}
\end{theorem}

\ni This is still a long way short of the conjecture that $r_k(N) < \pi(N)$ for $N$ sufficiently large. However, in addition to coming much closer to this bound than any previous arguments, Gowers succeeded in introducing methods of harmonic analysis to the problem for the first time since Roth. Since harmonic analysis (in the form of the circle method of Hardy and Littlewood) has been the most effective tool in tackling additive problems involving the primes, it seems fair to say that it was the work of Gowers which first gave us hope of tackling long progressions of primes. The ideas of Gowers will feature fairly substantially in this exposition, but in our paper \cite{green-tao-primes} much of what is done is more in the ergodic-theoretic spirit of Furstenberg and of more recent authors in that area such as Host--Kra \cite{host-kra} and Ziegler \cite{ziegler}. \vs

\ni To conclude this discussion of Szemer\'edi's theorem we mention a variant of it which is far more useful in practice. This applies to \textit{functions}\footnote{When discussing additive problems it is often convenient to work in the context of a finite abelian group $G$. For problems involving $\{1,\dots,N\}$ there are various technical tricks which allow one to work in $\Z/N'\Z$, for some $N' \approx N$. In this expository article we will not bother to distinguish between $\{1,\dots,N\}$ and $\Z/N\Z$. For examples of the technical trickery required here, see \cite[Definition 9.3]{green-tao-primes}, or the proof of Theorem 2.6 in \cite{gowers}.} $f : \Z/N\Z \rightarrow [0,1]$ rather than just to (characteristic functions of) sets. It also guarantees \textit{many} arithmetic progressions of length $k$. This version does, however, follow from the earlier formulation by some fairly straightforward averaging arguments due to Varnavides \cite{varnavides}.

\begin{proposition}[Szemer\'edi's theorem, II]\label{szem2}
Let $k \geq 3$ be an integer, and let $\delta \in (0,1]$ be a real number. Then there is a constant $c(k,\delta) > 0$ such that for any function $f : \Z/N\Z \rightarrow [0,1]$ with $\E f = \delta$ we have the bound\footnote{We use this very convenient conditional expectation notation repeatedly. $\E_{x \in A} f(x)$ is defined to equal $|A|^{-1}\sum_{x \in A} f(x)$.}
\begin{equation}\boxeq \E_{x,d \in \Z/N\Z} f(x) f(x+d) \dots f(x + (k-1)d) \geq c(k,\delta).\end{equation}
\end{proposition}

\ni We do not, in \cite{green-tao-primes}, prove any new bounds for $r_k(N)$. Our strategy is to prove a \textit{relative Szemer\'edi theorem}. To describe this we consider, for brevity of exposition, only the case $k = 4$. Consider the following table.\vs

\begin{center}
\begin{tabular}{|c||c|}
\hline
 Szemer\'edi & Relative Szemer\'edi\\
\hline
$\{1,\dots,N\}$ & ? \\
\hline
\begin{tabular}{c}
$A \subseteq \{1,\dots,N\}$ \\ $|A| \geq 0.0001 N$ \end{tabular} & \begin{tabular}{c}
$\mathcal{P}_N$ \\ = primes $\leq N$ \end{tabular} \\
\hline
\begin{tabular}{c} Szemer\'edi's theorem:\\ $A$ contains many $4$-term APs.\end{tabular} & \begin{tabular}{c} Green--Tao theorem: \\ $\mathcal{P}_N$ contains many $4$-term APs. \end{tabular}\\
\hline
\end{tabular}
\end{center}\vs
On the left-hand side of this table is Szemer\'edi's theorem for progressions of length $4$, stated as the result that a set $A \subseteq \{1,\dots,N\}$ of density $0.0001$ contains many $4$-term APs if $N$ is large enough. On the right is the result we wish to prove. Only one thing is missing: we must find an object to play the r\^ole of $\{1,\dots,N\}$. We might try to place the primes inside some larger set $\mathcal{P}'_N$ in such a way that $|\mathcal{P}_N| \geq 0.0001 |\mathcal{P}'_N|$, and hope to prove an analogue of Szemer\'edi's theorem for $\mathcal{P}'_N$. \vs

\ni A natural candidate for $\mathcal{P}'_N$ might be the set of \textit{almost primes}; perhaps, for example, we could take $\mathcal{P}'_N$ to be the set of integers in $\{1,\dots,N\}$ with at most 100 prime factors. This would be consistent with the intuition, coming form sieve theory, that almost primes are much easier to deal with than primes. It is relatively easy to show, for example, that there are long arithmetic progressions of almost primes \cite{grosswald}.\vs

\ni This idea does not quite work, but a variant of it does. Instead of a \textit{set} $\mathcal{P}'_N$ we instead consider what we call a \textit{measure}\footnote{Actually, $\nu$ is just a function but we use the term ``measure'' to distinguish it from other functions appearing in our work.} $\nu : \{1,\dots,N\} \rightarrow [0,\infty)$. Define the \textit{von Mangoldt function} $\Lambda$ by
\[ \Lambda(n) := \left\{ \begin{array}{ll}  \log p & \mbox{if $n = p^k$ is prime} \\ 0 & \mbox{otherwise}.\end{array}\right.\]
The function $\Lambda$ is a weighted version of the primes; note that the prime number theorem is equivalent to the fact that $\E_{1 \leq n \leq N} \Lambda(n) = 1 + o(1)$. Our measure $\nu$ will satisfy the following two properties.\vs

\begin{enumerate} \item ($\nu$ majorizes the primes) We have $\Lambda(n) \leq 10000\nu(n)$ for all $1 \leq n \leq N$.
\item  (primes sit inside $\nu$ with positive density) We have $\E_{1 \leq n \leq N} \nu(n) = 1 + o(1)$.
\end{enumerate}

\ni These two properties are very easy to satisfy, for example by taking $\nu = \Lambda$, or by taking $\nu$ to be a suitably normalised version of the almost primes. Remember, however, that we intend to prove a Szemer\'edi theorem relative to $\nu$. In order to do that it is reasonable to suppose that $\nu$ will need to meet more stringent conditions. The conditions we use in \cite{green-tao-primes} are called the \textit{linear forms condition} and the \textit{correlation condition}. We will not state them here in full generality, referring the reader to \cite[\S 3]{green-tao-primes} for full details. We remark, however, that verifying these conditions is of the same order of difficulty as obtaining asymptotics for, say,
\[ \sum_{n \leq N} \nu(n) \nu(n+2).\]
For this reason there is no chance that we could simply take $\nu = \Lambda$, since if we could do so we would have solved the twin prime conjecture.\vs

\ni We call a measure $\nu$ which satisfies the linear forms and correlation conditions \textit{pseudorandom}.\vs

\ni To succeed with the relative Szemer\'edi strategy, then, our aim is to find a pseudorandom measure $\nu$ for which conditions (i) and (ii) and the  are satisfied. Such a function\footnote{Actually, this is a lie. There is no pseudorandom measure which majorises the primes themselves. One must first use a device known as the $W$-trick to remove biases in the primes coming from their irregular distribution in residue classes to small moduli. This is discussed in \S \ref{sec3}.} comes to us, like the almost primes, from the idea of using a sieve to bound the primes. The particular sieve we had recourse to was the $\Lambda^2$-sieve of Selberg. Selberg's great idea was as follows.\vs

\ni Fix a parameter $R$, and let $\lambda = (\lambda_d)_{d = 1}^R$ be any sequence of real numbers with $\lambda_1 = 1$. Then the function

\[ \sigma_{\lambda}(n) := (\sum_{\substack{ d | n \\ d \leq R}} \lambda_d )^2\]
majorizes the primes greater than $R$. Indeed if $n > R$ is prime then the truncated divisor sum over $d | n$, $d \leq R$ contains just one term corresponding to $d = 1$. \vs

\ni Although this works for any sequence $\lambda$, some choices are much better than others. If one wishes to minimise
\[ \sum_{n \leq N} \sigma_{\lambda}(n)\]
then, provided that $R$ is a bit smaller than $\sqrt{N}$, one is faced with a minimisation problem involving a certain quadratic form in the $\lambda_d$s. The optimal weights $\lambda_d^{\mbox{\scriptsize SEL}}$, Selberg's weights, have a slightly complicated form, but roughly we have
\[ \lambda_d^{\mbox{\scriptsize SEL}} \approx \lambda_d^{\mbox{\scriptsize GY}} := \mu(d) \frac{\log(R/d)}{\log R},\]
where $\mu(d)$ is the M\"obius function. These weights were considered by Goldston and Y{\i}ld{\i}r{\i}m \cite{goldston-yildirim} in some of their work on small gaps between primes (and earlier, in other contexts, by others including Heath-Brown). It seems rather natural, then, to define a function $\nu$ by 
\[ \nu(n) := \left\{\begin{array}{ll} \log N & n \leq R \\ \displaystyle\frac{1}{\log R}\bigg( \sum_{\substack{d | n \\ d \leq R}}\lambda_d^{\mbox{\scriptsize GY}}\bigg)^2 & n > R.\end{array}\right.\]
The weight $1/\log R$ is chosen for normalisation purposes; if $R < N^{1/2 - \epsilon}$ for some $\epsilon > 0$ then we have $\E_{1 \leq n \leq N} \nu(n) = 1 + o(1)$. \vs

\ni One may more-or-less read out of the work of Goldston and Y{\i}ld{\i}r{\i}m a proof of properties (i) and (ii) above, as well as pseudorandomness, for this function $\nu$. One requires that $R < N^{c}$ where $c$ is sufficiently small. These verifications use the classical zero-free region for the $\zeta$-function and classical techniques of contour integration. \vs

\ni Goldston and Y{\i}ld{\i}r{\i}m's work was part of their long-term programme to prove that
\begin{equation}\label{eq1} \mbox{liminf}_{n \rightarrow \infty} \frac{p_{n+1} - p_n}{\log n} = 0,\end{equation}
where $p_n$ is the $n$th prime. We have recently learnt that this programme has been successful. Indeed together with J.~Pintz they have used weights coming from a higher-dimensional sieve in order to establish \eqref{eq1}. It is certain that without the earlier preprints of Goldston and Y{\i}ld{\i}r{\i}m our work would have developed much more slowly, at the very least.\vs

\ni Let us conclude this section by remarking that $\nu$ will not play a great r\^ole in the subsequent exposition. It plays a substantial r\^ole in \cite{green-tao-primes}, but in a relatively non-technical exposition like this it is often best to merely remark that the measure $\nu$ and the fact that it is pseudorandom is used all the time in proofs of the various statements that we will describe.

\section{progressions of length three and linear bias}
\label{sec3}

\ni Let $G$ be a finite abelian group with cardinality $N$. If $f_1,\dots,f_k : G \rightarrow \mathbb{C}$ are any functions we write
\[ T_k(f_1,\dots,f_k) := \E_{x,d \in G} f_1(x) f_2(x+d) \dots f_k(x + (k-1)d)\]
for the normalised count of $k$-term APs involving the $f_i$. When all the $f_i$ are equal to some function $f$, we write
\[ T_k(f) := T_k(f,\dots,f).\]
When $f$ is equal to $1_A$, the characteristic function of a set $A \subseteq G$, we write
\[ T_k(A) := T_k(1_A) = T_k(1_A,\dots,1_A).\]
This is simply the number of $k$-term arithmetic progressions in the set $A$, divided by $N^2$.\vs

\ni Let us begin with a discussion of 3-term arithmetic progressions and the trilinear form $T_3$. If $A \subseteq G$ is a set, then clearly $T_3(A)$ may vary between $0$ (when $A = \emptyset$) and $1$ (when $A = G$). If, however, one places some restriction on the cardinality of $A$ then the following question seems natural:

\begin{question}\label{ques3.1}
Let $\alpha \in (0,1)$, and suppose that $A \subseteq G$ is a set with cardinality $\alpha N$. What is $T_3(A)$?
\end{question}

\ni To think about this question, we consider some examples.\vs

\ni\textit{Example 1} (Random set). Select a set $A \subseteq G$ by picking each element $x \in G$ to lie in $A$ independently at random with probability $\alpha$. Then with high probability $|A| \approx \alpha N$. Also, if $d \neq 0$, the arithmetic progression $(x, x+d, x+ 2d)$ lies in $G$ with probability $\alpha^3$. Thus we expect that $T_3(A) \approx \alpha^3$, and indeed it can be shown using simple large deviation estimates that this is so with high probability.\vs

\ni Write $E_3(\alpha) := \alpha^3$ for the \textit{expected} normalised count of three-term progressions in the random set of Example 1. One might refine Question \ref{ques3.1} by asking:

\begin{question}\label{ques3.2}
Let $\alpha \in (0,1)$, and suppose that $A \subseteq G$ is a set with cardinality $\alpha N$. Is $T_3(A) \approx E_3(\alpha)$?
\end{question}

\ni It turns out that the answer to this question is ``no'', as the next example illustrates.\vs

\ni\textit{Example 2} (Highly structured set, I). Let $G = \Z/N\Z$, and consider the set $A = \{1,\dots, \lfloor\alpha N\rfloor\}$, an interval. It is not hard to check that if $\alpha < 1/2$ then $T_3(A) \approx \frac{1}{4}\alpha^2$, which is much bigger than $E_3(\alpha)$ for small $\alpha$.\vs

\ni These first two examples do not rule out a positive answer to the following question.

\begin{question}\label{ques3.3}
Let $\alpha \in (0,1)$, and suppose that $A \subseteq G$ is a set with cardinality $\alpha N$. Is $T_3(A) \geq E_3(\alpha)$?
\end{question}

\ni If this question did have an affirmative answer, the quest for progressions of length three in sets would be a fairly simple one (the primes would trivially contain many three-term progressions on density grounds alone, for example). Unfortunately, there are counterexamples.\vs

\ni\textit{Example 3} (Highly structured set, II). Let $G = \Z/N\Z$. Then there are sets $A \subseteq G$ with $|A| = \lfloor\alpha N\rfloor$ , yet with $T_3(A) \ll \alpha^{10000}$. We omit the details of the construction, remarking only that such sets can be constructed\footnote{Basically one considers a set $S \subseteq \Z^2$ formed as the product of a Behrend set in $\{1,\dots,M\}$ and the interval $\{1,\dots,L\}$, for suitable $M$ and $L$, and then one projects this set linearly to $\Z/N\Z$.} as unions of intervals of length $\gg_{\alpha} N$ in $\Z/N\Z$. \vs

\ni Our discussion so far seems to be rather negative, in that our only conclusion is that none of Questions \ref{ques3.1}, \ref{ques3.2} and \ref{ques3.3} have particularly satisfactory answers. Note, however, that the three examples we have mentioned are all consistent with the following dichotomy.

\begin{dichotomy}[Randomness vs Structure for $3$-term APs] \label{dic3}
Suppose that $A \subseteq G$ has size $\alpha N$. Then \textbf{either}
\begin{itemize}
\item $T_3(A) \approx E_3(\alpha)$ \textbf{or}
\item $A$ has \textit{structure}.
\end{itemize}
\end{dichotomy}

\ni It turns out that one may clarify, in quite a precise sense, what is meant by \textit{structure} in this context. The following proposition may be proved by fairly straightforward harmonic analysis. We use the Fourier transform on $G$, which is defined as follows. If $f : G \rightarrow \mathbb{C}$ is a function and $\gamma \in \widehat{G}$ a character (i.e. a homomorphism from $G$ to $\mathbb{C}^{\times}$), then
\[ f^{\wedge}(\gamma) := \E_{x \in G} f(x) \gamma(x).\]

\begin{proposition}[Too many/few 3APs implies linear bias]\label{3-dichotomy}
Let $\alpha, \eta \in (0,1)$. Then there is $c(\alpha,\eta) > 0$ with the following property. Suppose that $A \subseteq G$ is a set with $|A| = \alpha N$, and that 
\[ |T_3(A) - E_3(\alpha)| \geq \eta.\]
Then there is some character $\gamma \in \widehat{G}$ with the property that 
\begin{equation}\boxeq |(1_A - \alpha)^{\wedge}(\gamma)| \geq c(\alpha, \eta).\end{equation}
\end{proposition}

\ni Note that when $G = \Z/N\Z$ every character $\gamma$ has the form $\gamma(x) = e(rx/N)$. It is the occurrence of the linear function $x \mapsto rx/N$ here which gives us the name \textit{linear bias}.\vs

\ni It is an instructive exercise to compare this proposition with Examples 1 and 2 above. In Example 2, consider the character $\gamma(x) = e(x/N)$. If $\alpha$ is reasonably small then all the vectors $e(x/N)$, $x \in A$, have large positive real part and so when the sum
\[ (1_A - \alpha)^{\wedge}(\gamma) = \E_{x \in \Z/N\Z} \widehat{1}_A(x) e(x/N)\]
is formed there is very little cancellation, with the result that the sum is large.\vs

\ni In Example 1, by contrast, there is (with high probability) considerable cancellation in the sum for $(1_A - \alpha)^{\wedge}(\gamma)$ for every character $\gamma$.

\section{linear bias and the primes}

\ni What use is Dichotomy \ref{dic3} for thinking about the primes? One might hope to use Proposition \ref{3-dichotomy} in order to count 3-term APs in some set $A \subseteq G$ by showing that $A$ does not have linear bias. One would then know that $T_3(A) \approx E_3(\alpha)$, where $|A| = \alpha N$. \vs

\ni Let us imagine how this might work in the context of the primes. We have the following proposition\footnote{There are two ways of proving this proposition. One uses classical harmonic analysis. For pointers to such a proof, which would involve establishing an $L^p$-restriction theorem for $\nu$ for some $p \in (2,3)$, we refer the reader to \cite{green-tao}. This proof uses more facts about $\nu$ than mere pseudorandomness. Alternatively, the result may be deduced from Proposition \ref{3-dichotomy} by a transference principle using the machinery of \cite[\S 6--8]{green-tao-primes}. For details of this approach, which is far more amenable to generalisation, see \cite{green-tao-prime4aps}. Note that Proposition \ref{primes-3-dichotomy} does not feature in \cite{green-tao-primes} and is stated here for pedagogical reasons only.}, which is an analogue of Proposition \ref{3-dichotomy}. In this proposition\footnote{Recall that we are being very hazy in distinguishing between $\{1,\dots,N\}$ and $\Z/N\Z$.}, $\nu : \Z/N\Z \rightarrow [0,\infty)$ is the Goldston-Y{\i}ld{\i}r{\i}m measure constructed in \S \ref{sec2}. 

\begin{proposition}\label{primes-3-dichotomy} Let $\alpha,\eta \in (0,2]$. Then there is $c(\alpha,\eta) > 0$ with the following propety. Let $f :\Z/N\Z \rightarrow \mathbb{R}$ be a function with $\E f = \alpha$ and such that $0 \leq f(x) \leq 10000\nu(x)$ for all $x \in \Z/N\Z$, and suppose that 
\[ |T_3(f) - E_3(\alpha)| \geq \eta.\]
Then 
\begin{equation}\label{eq44} |\E_{x \in \Z/N\Z}(f(x) - \alpha)e(rx/N)| \geq c(\alpha,\eta)\end{equation} for some $r \in \Z/N\Z$.\endproof
\end{proposition}
\ni This proposition may be applied with $f = \Lambda$ and $\alpha = 1 + o(1)$. If we could rule out \eqref{eq44}, then we would know that $T_3(\Lambda) \approx E_3(1) = 1$, and would thus have an asymptotic for 3-term progressions of primes.\vs

\ni Sadly, \eqref{eq44} does hold. Indeed if $N$ is even and $r = N/2$ then, observing the most primes are odd, it is easy to confirm that
\[ \E_{x \in \Z/N\Z}(\Lambda(x) - 1) e(rx/N) = -1 + o(1).\]
That is, the primes \textit{do} have linear bias.\vs

\ni Fortunately, it is possible to modify the primes so that they have no linear bias using a device that we refer to as the $W$-trick. We have remarked that most primes are odd, and that as a result $\Lambda - 1$ has considerable linear bias. However, if one takes the odd primes 
\[ 3,5,7,11,13,17,19,\dots\]
and then rescales by the map $x \mapsto (x - 1)/2$, one obtains the set
\[ 1,2,3,5,6,8,9,\dots\]
which does not have substantial $\md{2}$ bias (this is a consequence of the fact that there are roughly the same number of primes congruent to $1$ and $3 \md{4}$). Furthermore, if one can find an arithmetic progression of length $k$ in this set of rescaled primes, one can certainly find such a progression in the primes themselves. Unfortunately this set of rescaled primes still has linear bias, because it contains only one element $\equiv 1 \md{3}$. However, a similar rescaling trick may be applied to remove this bias too, and so on. \vs

\ni Here, then, is the $W$-trick. Take a slowly growing function $w(N) \rightarrow \infty$, and set $W := \prod_{p < w(N)} p$. Define the rescaled von Mongoldt function $\widetilde{\Lambda}$ by

\[ \widetilde{\Lambda}(n) := \frac{\phi(W)}{W} \Lambda(Wn + 1).\]
The normalisation has been chosen so that $\E \widetilde{\Lambda} = 1 + o(1)$. $\widetilde{\Lambda}$ does not have substantial bias in any residue class to modulus $q < w(N)$, and so there is at least hope of applying a suitable analogue of Proposition \ref{primes-3-dichotomy} to it.\vs

\ni Now it is a straightforward matter to define a new pseudorandom measure  $\widetilde{\nu}$ which majorises $\widetilde{\Lambda}$. Specifically, we have

\begin{enumerate}
\item ($\widetilde{\nu}$ majorizes the modified primes) We have $\widetilde{\lambda}(n) \leq 10000\widetilde{\nu}(n)$ for all $1 \leq n \leq N$.
\item (modified primes sit inside $\widetilde{\nu}$ with positive density) We have $\E_{1 \leq n \leq N} \widetilde{\nu}(n) = 1 + o(1)$.\end{enumerate}

\ni The following modified version of Proposition \ref{primes-3-dichotomy} may be proved:

\begin{proposition}\label{primes-3-dichotomy-2} Let $\alpha,\eta \in (0,2]$. Then there is $c(\alpha,\eta) > 0$ with the following propety. Let $f :\Z/N\Z \rightarrow \mathbb{R}$ be a function with $\E f = \alpha$ and such that $0 \leq f(x) \leq 10000\widetilde{\nu}(x)$ for all $x \in \Z/N\Z$, and suppose that 
\[ |T_3(f) - E_3(\alpha)| \geq \eta.\]
Then 
\begin{equation}\label{eq44a} |\E_{x \in \Z/N\Z}(f(x) - \alpha)e(rx/N)| \geq c(\alpha,\eta)\end{equation} for some $r \in \Z/N\Z$. \endproof
\end{proposition}
\ni This may be applied with $f = \widetilde{\Lambda}$ and $\alpha = 1 + o(1)$. Now, however, condition \eqref{eq44a} does not so obviously hold. In fact, one has the estimate
\begin{equation}\label{eq44b} \sup_{r \in \Z/N\Z} |\E_{x \in \Z/N\Z}(\widetilde{\Lambda}(x) - 1)e(rx/N)| = o(1).\end{equation}
To prove this requires more than simply the good distribution of $\widetilde{\Lambda}$ in residue classes to small moduli. It is, however, a fairly standard consequence of the Hardy-Littlewood circle method as applied to primes by Vinogradov. In fact, the whole theme of linear bias in the context of additive questions involving primes may be traced back to Hardy and Littlewood.\vs

\ni Proposition \ref{primes-3-dichotomy-2} and \eqref{eq44b} imply that $T_3(\widetilde{\Lambda}) \approx E_3(1) = 1$. Thus there are infinitely many three-term progressions in the modified ($W$-tricked) primes, and hence also in the primes themselves\footnote{In fact, this analysis does not have to be pushed much further to get a proof of Conjecture \ref{hlconj} for $k=3$, that is to say an asymptotic for 3-term progressions of primes. One simply counts progressions $x,x+d,x+2d$ by splitting into residue classes $x \equiv b \md{W}$, $d \equiv b' \md{W}$ and using a simple variant of Proposition \ref{primes-3-dichotomy-2}.}.

\section{progressions of length four and quadratic bias}

\ni We return now to the discussion of \S \ref{sec3}. There we were interested in counting 3-term arithmetic progressions in a set $A \subseteq G$ with cardinality $\alpha N$. In this section our interest will be in $4$-term progressions.\vs

\ni Suppose then that $A \subseteq G$ is a set, and recall that
\[ T_4(A) := \E_{x,d  \in G} 1_A(x)1_A(x+d)1_A(x+2d)1_A(x + 3d)\] is the normalised count of four-term arithmetic progressions in $A$. One may, of course, ask the analogue of Question \ref{ques3.1}:

\begin{question}\label{ques4.1}
Let $\alpha \in (0,1)$, and suppose that $A \subseteq G$ is a set with cardinality $\alpha N$. What is $T_4(A)$?
\end{question}

\ni Examples 1,2 and 3 make perfect sense here, and we see once again that there is no immediately satisfactory answer to Question \ref{ques4.1}. With high probability the random set of Example 1 has about $E_4(\alpha) := \alpha^4$ four-term APs, but there are structured sets with with substantially more or less than this number of APs. As in \S \ref{sec3}, these examples are consistent with a dichotomy of the following type:

\begin{dichotomy}[Randomness vs Structure for $4$-term APs] \label{dic4}
Suppose that $A \subseteq G$ has size $\alpha N$. Then \textbf{either}
\begin{itemize}
\item $T_4(A) \approx E_4(\alpha)$ \textbf{or}
\item $A$ has \textit{structure}.
\end{itemize}
\end{dichotomy}

Taking into account the three examples we have so far, it is quite possible that this dichotomy takes \textit{exactly} the form of that for 3-term APs. That is to say ``$A$ has structure'' could just mean that $A$ has linear bias:

\begin{question}\label{4-linear}
Let $\alpha, \eta \in (0,1)$.  Suppose that $A \subseteq G$ is a set with $|A| = \alpha N$, and that 
\[ |T_4(A) - E_4(\alpha)| \geq \eta.\]
Must there exist some $c = c(\alpha, \eta) > 0$ and some character $\gamma \in \widehat{G}$ with the property that 
\[ |(1_A - \alpha)^{\wedge}(\gamma)| \geq c(\alpha, \eta)?\]
\end{question}

\ni That the answer to this question is no, together with the nature of the counterexample, is one of the key themes of our whole work. This phenomenon was discovered, in the context of ergodic theory, by Furstenberg and Weiss \cite{furst-weiss} and then again, in the discrete setting, by Gowers \cite{gowers}. \vs

\ni\textit{Example 4} (Quadratically structured set). Define $A \subseteq \Z/N\Z$ to be the set of all $x$ such that $x^2 \in [-\alpha N/2, \alpha N/2]$. It is not hard to check using estimates for Gauss sums that $|A| \approx \alpha N$, and also that 
\[ \sup_{r \in \Z/N\Z} |\E(1_A(x) - \alpha)e(rx/N)| = o(1),\] that is to say $A$ does not have linear bias.
(In fact, the largest Fourier coefficient of $1_A - \alpha$ is just $N^{-1/2 + \epsilon}$.) Note, however, the relation
\[ x^2 - 3(x+d)^2 + 3(x + 2d)^2 + (x + 3d)^2 = 0,\]
valid for arbitrary $x,d \in \Z/N\Z$. This means that if $x, x+ d, x+2d \in A$ then automatically we have
\[ (x + 3d)^2 \in [-7\alpha N/2, 7\alpha N/2].\]
It seems, then, that if we know that $x, x+d$ and $x + 2d$ lie in $A$ there is a very high chance that $x + 3d$ also lies in $A$. This observation may be made rigorous, and it does indeed transpire that $T_4(A) \geq c\alpha^3$.\vs

\ni How can one rescue the randomness-structure dichotomy in the light of this example? Rather remarkably, ``quadratic'' examples like Example 4 are the \textit{only} obstructions to having $T_4(A) \approx E_4(\alpha)$. There is an analogue of Proposition \ref{3-dichotomy} in which characters $\gamma$ are replaced by ``quadratic'' objects\footnote{The proof of this proposition is long and difficult and may be found in \cite{green-tao-gowersu3}. It is heavily based on the arguments of Gowers \cite{gowers4,gowers}. This proposition has no place in \cite{green-tao-primes}, and it is once again included for pedagogical reasons only. It played an important r\^ole in the development of our ideas.}. 

\begin{proposition}[Too many/few 4APs implies quadratic bias]\label{4-dichotomy}
Let $\alpha, \eta \in (0,1)$. Then there is $c(\alpha,\eta) > 0$ with the following property. Suppose that $A \subseteq G$ is a set with $|A| = \alpha N$, and that 
\[ |T_4(A) - E_4(\alpha)| \geq \eta.\]
Then there is some quadratic object $q \in \mathcal{Q}(\kappa)$, where $\kappa \geq \kappa_0(\alpha,\eta)$, with the property that 
\begin{equation}\boxeq |\E_{x \in G}(1_A(x) - \alpha)q(x)| \geq c(\alpha, \eta).\end{equation}
\end{proposition}

\ni We have not, of course, said what we mean by the set of \textit{quadratic objects} $\mathcal{Q}(\kappa)$. To give the exact definition, even for $G = \Z/N\Z$, would take us some time, and we refer to \cite{green-tao-gowersu3} for a full discussion. In the light of Example 4, the reader will not be surprised to hear that quadratic exponentials such as $q(x) = e(x^2/N)$ are members of $\mathcal{Q}$. However, $\mathcal{Q}(\kappa)$ also contains rather more obscure objects\footnote{We are thinking of these as defined on $\{1,\dots,N\}$ rather than $\Z/N\Z$.} such as 
\[ q(x) = e(x\sqrt{2} \{ x \sqrt{3}\})\]
and
\[ q(x) = e(x\sqrt{2} \{x \sqrt{3}\} + x \sqrt{5} \{x \sqrt{7}\} + x\sqrt{11}),\]
where $\{x\}$ denotes fractional part.  The parameter $\kappa$ governs the complexity of the expressions which are allowed: smaller values of $\kappa$ correspond to more complicated expressions. The need to involve these ``generalised'' quadratics in addition to ``genuine'' quadratics such as $e(x^2/N)$ was first appreciated by Furstenberg and Weiss in the ergodic theory context, and the matter also arose in the work of Gowers. 

\section{quadratic bias and the primes}

\ni It is possible to prove\footnote{As with Proposition \ref{primes-3-dichotomy}, this proposition does not appear in \cite{green-tao-primes}, though it motivated our work and a variant of it is used in our later work \cite{green-tao-prime4aps}. Once again there are two proofs. One is based on a combination of harmonic analysis and the work of Gowers, is difficult, and requires more facts about $\nu$ than mere pseudorandomness. This was our original argument. It is also possible to proceed by a transference principle, deducing the result from Proposition \ref{4-dichotomy} using the machinery of \cite[\S 6--8]{green-tao-primes}. See \cite{green-tao-prime4aps} for more details.} a version of Proposition \ref{4-dichotomy} which might be applied to primes. The analogue of Proposition \ref{primes-3-dichotomy} is true but not useful, for the same reason as before: the primes exhibit significant bias in residue classes to small moduli. As before, this bias may be removed using the $W$-trick. 

\begin{proposition}\label{primes-4-dichotomy-2} Let $\alpha,\eta \in (0,2]$. Then there are $c(\alpha,\eta)$ and $\kappa_0(\alpha,\eta) > 0$ with the following propety. Let $f :\Z/N\Z \rightarrow \mathbb{R}$ be a function with $\E f = \alpha$ and such that $0 \leq f(x) \leq 10000\widetilde{\nu}(x)$ for all $x \in \Z/N\Z$, and suppose that 
\[ |T_4(f) - E_4(\alpha)| \geq \eta.\]
Then we have
\begin{equation}\label{eq44d} |\E_{x \in \Z/N\Z}(f(x) - \alpha)q(x)| \geq c(\alpha,\eta)\end{equation} for some quadratic object $q \in \mathcal{Q}(\kappa)$ with $\kappa \geq \kappa_0(\alpha,\eta)$.\endproof
\end{proposition}

\ni One is interested, of course, in applying this with $f = \widetilde{\Lambda}$. If we could verify that \eqref{eq44d} does not hold, that is to say the primes do not have quadratic bias, then it would follow that $T_4(\widetilde{\Lambda}) \approx E_4(1) = 1$. This means that the modified ($W$-tricked) primes have many 4-term progressions, and hence so do the primes themselves\footnote{in fact, just as for progressions of length 3, this allows one to obtain a proof of Conjecture \ref{hlconj} for $k = 4$, that is to say an asymptotic for prime progressions of length 4. See \cite{green-tao-prime4aps}.}.\vs

\ni One wishes to show, then, that for fixed $\kappa$ one has 
\begin{equation}\label{lam-quad} \sup_{q \in \mathcal{Q}(\kappa)}|\E_{x \in \Z/N\Z} (\widetilde{\Lambda}(x) - 1)q(x) | = o(1).\end{equation}
Such a result is certainly not a consequence of the classical Hardy-Littlewood circle method\footnote{Though reasonably straightforward extensions of the circle method do permit one to handle genuine quadratic phases such as $q(x) = e(x^2\sqrt{2})$.}. Generalised quadratic phases such as $q(x) = e(x\sqrt{2} \{x \sqrt{3}\})$ are particularly troublesome. Although we do now have a proof of \eqref{lam-quad}, it is very long and complicated. See \cite{green-tao-u3mobius} for details. \vs

\ni In the next section we explain how our original paper \cite{green-tao-primes} managed to avoid the need to prove \eqref{lam-quad}.

\section{quotienting out the bias - the energy increment argument}\label{sec7}

\ni Our paper \cite{green-tao-primes} failed to rule out the possibility that $\widetilde{\Lambda} - 1$ correlates with some quadratic function $q \in \mathcal{Q}(\kappa)$. For that reason we did not obtain a proof of Conjecture \ref{hlconj}, getting instead the weaker statement of Theorem \ref{gtthm}. In this section\footnote{The exposition in this section is rather looser than in other sections. To make the argument rigorous, one must introduce various technical devices, such as the exceptional sets which feature in \cite[\S 7,8]{green-tao-primes}. We are also being rather vague about the meaning of terms such as ``correlate'', and the parameter $\kappa$ involved in the definition of quadratic object. Note also that the argument of \cite{green-tao-primes} uses \textit{soft} quadratic objects rather than the genuine ones which we are discussing here for expositional purposes. See \S \ref{sec8} for a brief discussion of these.} we outline the \textit{energy increment} argument of \cite{green-tao-primes}, which allowed us to deal with the possibility that $\widetilde{\Lambda} - 1$ does correlate with a quadratic. \vs

\ni We begin by writing
\begin{equation}\label{step0} \widetilde{\Lambda} := 1 + f_0.\end{equation}
Proposition \ref{primes-4-dichotomy-2} tells us that $T_4(\widetilde{\Lambda}) \approx 1$, unless $f_0$ correlates with some quadratic $q_0 \in \mathcal{Q}$. Suppose, then, that 
\[ |\E_{x \in \Z/N\Z} f_0(x) q_0(x) | \geq \eta.\]
Then we revise the decomposition \eqref{step0} to
\begin{equation}\label{step1} \widetilde{\Lambda} := F_1 + f_1,\end{equation}
where $F_1$ is a function defined using $q_0$. In fact, $F_1$ is basically the average of $\widetilde{\Lambda}$ over approximate level sets of $q_0$. That is, one picks an appropriate scale\footnote{As we remarked, the actual situation is more complicated. There is an averaging over possible decompositions of $[0,1]$ into intervals of length $\epsilon$, to ensure that the level sets look pleasant. There is also a need to consider exceptional sets, which unfortunately makes the argument look rather messy.} $\epsilon = 1/J$, and then defines 
\[ F_1 := \E( \widetilde{\Lambda} | \mathcal{B}_0),\]
where $\mathcal{B}_0$ is the $\sigma$-algebra generated by the sets $x : q_0(x) \in [j/J, (j+1)/J)$.\vs

\ni A variant of Proposition \ref{primes-4-dichotomy-2} implies a new dichotomy: either
$T_4(\widetilde{\Lambda}) \approx T_4(F_1)$, or else $f_1$ correlates with some quadratic $q_1 \in \mathcal{Q}$. Suppose then that
\[ |\E_{x \in \Z/N\Z} f_1(x) q_1(x)| \geq \eta.\]
We then further revise the decomposition \eqref{step1} to
\[ \widetilde{\Lambda} := F_2 + f_2,\]
where now 
\[ F_2 := \E(\widetilde{\Lambda} | \mathcal{B}_0 \wedge \mathcal{B}_1),\] the $\sigma$-algebra being defined by the joint level sets of $q_0$ and $q_1$.  \vs

\ni We repeat this process. It turns out that the algorithm stops in a finite number $s$ of steps, bounded in terms of $\eta$. the reason for this is that each new assumption
\[ |\E_{x \in \Z/N\Z} f_j(x) q_j(x) | \geq \eta\]
implies an increased lower bound for the \textit{energy} of $\widetilde{\Lambda}$ relative to the $\sigma$-algebra $\mathcal{B}_0 \wedge \dots \wedge \mathcal{B}_{j-1}$, that is to say the quantity
\[ E_j := \Vert \E(\widetilde{\Lambda} | \mathcal{B}_0 \wedge \dots \wedge \mathcal{B}_{j-1})\Vert_2.\]
The fact that $\widetilde{\Lambda}$ is dominated by $\widetilde{\nu}$ does, however, provide a universal bound for the energy, by dint of the evident inequality
\[ E_j \leq 10000\Vert \E(\widetilde{\nu} | \mathcal{B}_0 \wedge \dots \wedge \mathcal{B}_{j-1})\Vert_2.\]
The pseudorandomness of $\widetilde{\nu}$ allows one\footnote{This deduction uses the machinery of the Gowers $U^3$-norm, which we do not discuss in this survey. See \cite[\S 6]{green-tao-primes} for a full discussion. Of specific relevance is the fact that $\Vert \widetilde{\nu} \Vert_{U^3} = o(1)$, which is a consequence of the pseudorandomness of $\widetilde{\nu}$.} to bound the right-hand side here by $O(1)$. \vs

\ni At termination, then, we have a decomposition
\[ \widetilde{\Lambda} = F_s + f_s,\]
where 
\begin{equation}\label{eq478} \sup_{q \in \mathcal{Q}} |\E_{x \in \Z/N\Z} f_s(x) q(x) | < \eta,\end{equation}
and $F_s$ is defined by
\begin{equation}\label{eq480} F_s := \E(\widetilde{\Lambda} | \mathcal{B}_0 \wedge \mathcal{B}_1 \wedge \dots \wedge \mathcal{B}_{s-1}).\end{equation}
A variant of Proposition \ref{primes-4-dichotomy-2} implies, together with \eqref{eq478}, that 
\begin{equation}\label{eq479}
T_4(\widetilde{\Lambda}) \approx T_4(F_s).
\end{equation}
What can be said about $T_4(F_s)$? Let us note two things about the function $F_s$. First of all the definition \eqref{eq480} implies that
\begin{equation}\label{property1} \E F_s = \E \widetilde{\Lambda} = 1 + o(1).\end{equation}
Secondly, $F_s$ is not too large pointwise; this is again an artifact of $\widetilde{\Lambda}$ being dominated by $\widetilde{\nu}$. We have, of course,
\[ \Vert F_s \Vert_{\infty} = \Vert \E(\widetilde{\Lambda} | \mathcal{B}_0 \wedge \mathcal{B}_1 \wedge \dots \wedge \mathcal{B}_{s-1}) \Vert_{\infty} \leq 10000\Vert \E(\widetilde{\nu} | \mathcal{B}_0 \wedge \mathcal{B}_1 \wedge \dots \wedge \mathcal{B}_{s-1}) \Vert_{\infty}.\]
The pseudorandomness of $\widetilde{\nu}$ can again be used\footnote{Again, the machinery of the Gowers $U^3$-norm is used.} to show that the right-hand side here is $10000 + o(1)$; that is,
\begin{equation}\label{property2}
\Vert F_s \Vert_{\infty} \leq 10000 + o(1).
\end{equation}
\ni The two properties \eqref{property1} and \eqref{property2} together mean that $F_s$ behaves rather like the characteristic function of a subset of $\Z/N\Z$ with density at least $1/10000$. This suggests the use of Szemer\'edi's theorem to bound $T_4(F_s)$ below. The formulation of that theorem given in Proposition \ref{szem2} applies to exactly this situation, and it tells us that 
\[ T_4(F_s) > c\]
for some absolute constant $c > 0$.  Together with \eqref{eq479} this implies a similar lower bound for $T_4(\widetilde{\Lambda})$, which means that there are infinitely many 4-term arithmetic progressions of primes.\vs

\ni Let us conclude this section with an overview of what it is we have proved. The only facts about $\widetilde{\Lambda}$ that we used were that it is dominated pointwise by $10000\widetilde{\nu}$, and that $\E \widetilde{\Lambda}$ is not too small. The argument sketched above applies equally well in the general context of functions with these properties, and in the context of an arbitrary pseudorandom measure (not just the Goldston-Y{\i}ld{\i}r{\i}m measure). 

\begin{proposition}[Relative Szemer\'edi Theorem]\label{prop7} Let $\delta \in (0,1]$ be a real number and let $\nu$ be a psuedorandom measure. Then there is a constant $c'(4,\delta) > 0$ with the following property.
Suppose that $f : \Z/N\Z \rightarrow \mathbb{R}$ is a function such that $0 \leq f(x) \leq \nu(x)$ pointwise, and for which $\E f \geq \delta$. Then we have the estimate
\begin{equation}\boxeq T_4(f) \geq c'(4,\delta).\end{equation}
\end{proposition}
\ni In \cite{green-tao-primes} we prove the same\footnote{Note, however, that the definition of pseudorandom measure is strongly dependent on $k$.} theorem for progressions of any length $k \geq 3$. \vs

\ni Proposition \ref{prop7} captures the spirit of our argument quite well. We first deal with arithmetic progressions in a rather general context. Only upon completion of that study do we concern ourselves with the primes, and this is simply a matter of constructing an appropriate pseudorandom measure. Note also that Szemer\'edi's theorem is used as a ``black box''. We do not need to understand the proof of it, or to have good bounds for it.\vs

\ni Observe that one consequence of Proposition \ref{prop7} is a \textit{Szemer\'edi theorem relative to the primes}: any subset of the primes with positive relative density contains progressions of arbitrary length. Applying this to the set of primes congruent to $1 \md{4}$, we see that there are arbitrarily long progressions of numbers which are sums of two squares.

\section{soft obstructions}\label{sec8}

\ni Readers familiar with \cite{green-tao-primes} may have been confused by our exposition thus far, since ``quadratic objects'' play essentially no r\^ole in that paper. The purpose of this brief section is to explain why this is so, and to provide a bridge between this survey and our paper. Further details and discussion may be found in \cite[\S 6]{green-tao-primes}.\vs

\ni Let us start by recalling \S \ref{sec3}, where a set of ``obstructions'' to a set $A \subseteq G$ having roughly $E_3(\alpha)$ three-term APs was obtained. This was just the collection of characters $\gamma \in \widehat{G}$, and we used the term \textit{linear bias} to describe correlation with one of these characters.\vs

\ni Let $f : G \rightarrow \mathbb{C}$ be a function with $\Vert f \Vert_{\infty} \leq 1$. Now we observe the formula
\[ \E_{a,b \in G} \overline{f(x+a)f(x+b)}f(x+a+b) = \sum_{\gamma \in \widehat{G}} |\widehat{f}(\gamma)|^2 \overline{\widehat{f}(\gamma)} \overline{\gamma(x)},\]
which may be verified by straightforward harmonic analysis on $G$. Coupled with the fact that 
\[ \sum_{\gamma \in \widehat{G}} |\widehat{f}(\gamma)|^2 \leq 1,\]
a consequence of Parseval's identity, this means that the\footnote{The subscript 2 refers to the Gowers $U^2$-norm, which is relevant to the study of progressions of length 3.} ``dual function''
\[ \mathcal{D}_2f := \E_{a,b \in G} \overline{f(x+a)f(x+b)}f(x+a+b)\]
can be approximated by the weighted sum of a few characters. Every character is actually equal to a dual function; indeed we clearly have $\mathcal{D}_2(\gamma) = \overline{\gamma}$. \vs

\ni We think of the dual functions $\mathcal{D}_2(f)$ as \textit{soft linear obstructions}. They may be used in the iterative argument of \S \ref{sec7} in place of the genuinely linear functions, after one has established certain algebraic closure properties of these functions (see \cite[Proposition 6.2]{green-tao-primes})\vs

\ni The great advantage of these soft obstructions is that it is reasonably obvious how they should be generalised to give objects appropriate for the study of longer arithmetic progressions. We define
\[ \mathcal{D}_3(f) := \E_{a,b,c} \overline{f(x+a)f(x+b)f(x+c)}f(x+a+b)f(x+a+c)f(x+b+c)\overline{f(x+a+b+c)}.\]
This is a kind of sum of $f$ over parallelepipeds (minus one vertex), whereas $\mathcal{D}_2(f)$ was a sum over parallelograms (minus one vertex). This we think of as a \textit{soft quadratic obstruction}. Gone are the complications of having to deal with explicit generalised quadratic functions which, rest assured, only become worse when one deals with progressions of length 5 and longer.\vs

\ni The idea of using these soft obstructions came from the ergodic-theory work of Host and Kra \cite{host-kra}, where very similar objects are involved.\vs

\ni We conclude by emphasising that soft obstructions lead to relatively soft results, such as Theorem \ref{gtthm}. To get a proof of Conjecture \ref{hlconj} it will be necessary to return to generalised quadratic functions and their higher-order analogues.

\section{acknowledgements} \ni I would like to thank James Cranch for reading the manuscript and advice on using Mathematica, and Terry Tao for several helpful comments.

\end{document}